\documentstyle[12pt]{article}
\newlength{\abstractwidth}
\renewcommand{\thefootnote}{\fnsymbol{footnote}}
\renewcommand{\thanks}[1]{\footnote{#1}} 
\newcommand{\starttext}{
\setcounter{footnote}{0}
\renewcommand{\thefootnote}{\arabic{footnote}}}

\newcommand{\be}{\begin{equation}}
\newcommand{\bea}{\begin{eqnarray}}
\newcommand{\eea}{\end{eqnarray}}
\newcommand{\ee}{\end{equation}}

\def\ba{\begin{eqnarray}}
\def\ea{\end{eqnarray}}

\def\v{\vskip .1in}

\def\al{\alpha}
\def\b{\beta}

\def\d{\delta}
\def\e{\epsilon}
\def\g{\gamma}
\def\l{\lambda}
\def\m{\mu}
\def\n{\nu}
\def\o{\omega}

\def\r{\rho}
\def\si{\sigma}
\def\t{\theta}

\def\O{\Omega}
\def\T{\Theta}

\def\cH{{\cal H}}

\def\cL{{\cal L}}

\def\cS{{\cal S}}

\def\cX{{\cal X}}

\def\Z{{\bf Z}}

\def\R{{\bf R}}
\def\C{{\bf C}}
\def\P{{\bf P}}
\def\K{{\rm K\"ahler }}

\def\Aut{{\rm Aut}}
\def\End{{\rm End}}
\def\osc{{\rm osc}}

\def\ti\tilde

\def\u{\underline}

\def\pl{\partial}

\def\i{\infty}
\def\I{\int}

\def\s{\sum}

\def\ddb{\partial\bar\partial}
\def\sub{\subseteq}
\def\ra{\rightarrow}

\def\L{\Lambda}

\def\us{{\underline s}}

\def\tr{{\rm tr}}
\def\det{{\rm det}}

\def\ti{\tilde}

\def\Tr{{\rm Tr}}

\def\o{\omega}

\def\[{{\bf [}}
\def\]{{\bf ]}}

\def\pl{\partial}

\v

\begin{document} \starttext \baselineskip=15pt
\setcounter{footnote}{0} \newtheorem{theorem}{Theorem}
\newtheorem{lemma}{Lemma} \newtheorem{corollary}{Corollary}
\newtheorem{definition}{Definition} \begin{center} {\Large \bf TEST
CONFIGURATIONS FOR K-STABILITY AND GEODESIC RAYS}
\footnote{Research
supported in part by National Science Foundation grants DMS-02-45371
and DMS-05-14003}
\\
\bigskip

{\large D.H. Phong$^*$ and
Jacob Sturm$^{\dagger}$} \\

\bigskip

$^*$ Department of Mathematics\\
Columbia University, New York, NY 10027\\

\v

$^{\dagger}$ Department of Mathematics \\
Rutgers University, Newark, NJ 07102

\end{center}
\v
\v\v\v \begin{abstract}

{\small Let $X$ be a compact complex manifold, $L\ra X$ an ample
line bundle over $X$, and  $\cH$ the space of all positively curved
metrics on $L$. We show that a pair $(h_0,T)$ consisting of a point
$h_0\in\cH$ and a test configuration $T=(\cL\ra \cX\ra \C)$,
canonically determines a weak geodesic ray $R(h_0,T)$ in $\cH$ which
emanates from $h_0$. Thus a test configuration behaves like a vector
field on the space of \K potentials $\cH$. We prove that $R$ is
non-trivial if the $\C^\times$ action on $X_0$, the central fiber of $\cal X$,
is non-trivial. The ray $R$ is obtained  as limit of
smooth geodesic rays $R_k\sub\cH_k$, where $\cH_k\sub\cH$ is the
subspace of Bergman metrics.}

\end{abstract}

\bigskip
\baselineskip=15pt
\setcounter{equation}{0}
\setcounter{footnote}{0}

\section{Introduction}
\setcounter{equation}{0}

Let $X$ be a compact complex manifold.
According to a basic conjecture of Yau \cite{Y87}, the existence
of canonical metrics on $X$ should be
equivalent to a stability condition
in the sense of geometric invariant theory. A version
of this conjecture, due to
Tian \cite{T97} and Donaldson \cite{D02}, says that if
$L\ra X$ is an ample line bundle, then $X$ has a
metric of constant scalar curvature in $c_1(L)$
if and only if the pair $(X,L)$
is K-stable, that is, if and only if the Futaki invariant $F(T)$ is
negative for each
non-trivial test configuration $T$. In particular,
$F(T)<0$ for all such $T$ should imply that the K-energy $\n: \cH \ra \R$
is bounded below, where
$\cH$  is the space of all positively curved metrics on
$L$.

\v

Now it is well known that the K-energy is convex along  geodesics of
$\cH$ (Donaldson \cite{D99}). Thus, if $h_0\in\cH$ and if
$R:(-\i,0]\ra \cH$ is a smooth geodesic ray emanating from $h_0$,
then the restriction of $\n$ to $R$ is a smooth convex function
$\n_R: (-\i,0]\ra\R$ and hence $\lim_{t\ra-\i} \dot\n_R = a(R)$ is
well defined (here $\dot \n_R$ is the time derivative of the K-energy).
In particular, if $a(R)<0$, then $\n$ is bounded below on
the ray $R$.

\v
We are thus led to the following plan for relating K-stability to lower
bounds for the K-energy: Given a non-trivial test
configuration $T=(\cal L\ra\cX\ra\C)$ and a point $h_0\in\cH$,

\medskip
A) Associate to $(h_0,T)$ a  canonical non-trivial geodesic ray $R(T,h_0)$
emanating from $h_0$.
\v
B) Prove that
$\lim_{t\to-\i}\dot\n_R=F(T)+d(T)$.

\v
where $d(T)\geq 0$ has the property: $d(T)=0$ if $X_0$, the central fiber of
$\cX$, has no multiplicity and $F(T)<0$ implies $F(T) + d(T)< 0$.
If this plan could be implemented, then $F(T)<0$ for a single test
configuration $T$ would
imply that $\n$ is bounded
below on the ray $R(T,h_0)$.
And the K-stability of $(X,L)$ would imply that
$\n$ is bounded below on all the rays $R(T,h_0)$ emanating from $h_0$.

\v

In this paper, we take a step in the direction of the
plan outlined above.: For step A), we start with an arbitrary test configuration $T$ and
an arbitrary point $h_0\in\cH$. We associate to this data a weak
geodesic $R(h_0,T)$ which is upper semi-continuous
(but may not be smooth). If the $\C^\times$
action on $X_0$
is non-trivial (in particular, if $F(T)\not=0$),
then we show that $R(h_0,T)$ is a non-trivial geodesic.
\v

We also provide evidence for step B): The ray  $R(h_0,T)$ is constructed as a limit of
Bergman geodesic rays $h(t;k)$.
Under certain geometric conditions (which
are necessary for our proofs, but we expect can be removed) we observe that the
limit of the K-energy time derivative along $h(t;k)=h_0e^{-\phi(t;k)}$
converges to the Futaki invariant $F(T)$
as $k\ra \i$ if $X_0$ is multiplicity free.
\v
After raising $\cL$ and $L$ to sufficiently high powers, we may assume
that $L$ is very ample, that $H^0(X,L)$ generates $\oplus_{k=0}^\i H^0(X,L^k)$,
and that $\cL$ has exponent one (note that raising the power of the
line bundle will just amount to a reparametrization of the geodesic).
These assumptions will be made throughout this
paper.
\v
Our main results are Theorem 1 and Theorem 2 below
(with relevant notation provided in \S 4):
\v

\begin{theorem}\label{main}
Let $L\ra X$ be a very ample line bundle, $h_0$ a positively curved metric on
$L$, and $T$ a test configuration
for $(X,L)$. Let

\be\label{geo}
\phi_t\ = \ \lim_{k\to\i}\big(\sup_{l\geq k} [\phi(t;l)]\big)^*
\ee
Then $h(t)=h_0e^{-\phi_t}$ is a weak geodesic ray emanating from $h_0$.
Here we make use of the notation $u^*(\zeta_0)=\lim_{\e\to 0}
\sup_{|\zeta-\zeta_0|<\e}u(\zeta)$
for any locally bounded
$u:X\times (-\i,0]\ra \R$.

\end{theorem}

\begin{theorem}
Assume that the action of $\C^\times$ on $X_0$ is non-trivial. Then
the weak geodesic defined by $\phi_t$ in Theorem 1 is non-trivial.
\end{theorem}

\v
We note that the $\C^\times$ action on $X_0$
is non-trivial
if the Futaki invariant $F(T)$ of the
test configuration $T$ does not vanish.

\v\v
We also note the following result:

\begin{theorem}
Assume that the test configuration can be equivariantly imbedded in a proper family
${\cal X}\ra B$, where $\cal X$ and $B$  are smooth compact manifolds with the property that
the Chern class map $Pic(B)\ra H^2(B,{\bf Z})$ is injective and $X_0$ is multiplicity free. Then,
for each $k>0$,
\be
\label{double limit}\lim_{t\to-\i}\dot \n_k =F(T)
\ee
Here $\n_k$ is the restriction of $\n$
to the Bergman geodesic $h(t;k)$.
\end{theorem}
{\bf Remark:} Theorem 1 holds in a wider context than that stated above - our
proofs show that one can associate a weak geodesic ray to an arbitrary traceless
hermitian matrix $A\in gl(H^0(X,L))$ (thus the eigenvalues of $A$ are real
numbers and not restricted to lie in $\Z$).

\v

To define what is meant by a weak geodesic, we start by
recalling that
$\cH$ is an infinite dimensional symmetric space with respect
to its natural Riemannian structure (see
Mabuchi \cite{M87}, Semmes \cite{S92} and Donaldson \cite{D99}).
Furthermore, the geodesic equation for $h_0e^{-\phi_t}$ is equivalent to
the degenerate Monge-Amp\`ere equation
\be\label{MA}
\O^{n+1}=0\quad {\rm on}\quad X\times A
\ee
where $A\sub\C$ is an annulus (in the case of a geodesic
segment) or a punctured disk (in the case of a geodesic ray).
Here $\O=\Omega_\phi$ is the smooth $(1,1)$-form
on $X\times A$ determined by: $\O=\O_0+{\sqrt{-1}\over 2}\ddb\Phi$ where
$\O_0=p_1^*\o_0$, $\o_0$ is the curvature of $h_0$,
$p_1(x,w)=x$, $\Phi(x,w)=\phi_t(x)$, and $t=\log |w|$.
A weak geodesic $\phi_t$
is one for which $\O_{\phi}$ is a plurisubharmonic solution
to (\ref{MA}) in the sense of pluripotential theory \cite{BT76}.
\v

The problem of constructing geodesic rays from test
configurations has been considered previously by
Arezzo-Tian \cite{AT}. They show that,
if the central fiber of the test configuration $T$ is smooth,
then one can use the Cauchy-Kowalevska theorem to find a local
analytic solution near infinity to the geodesic equation, and in this way, they
construct a geodesic ray $R(T)$ in $\cH$.
In fact, they construct a
family of  rays $R_j(T)$ where $j$ ranges over certain free
parameters which determine the power series coefficients.
These rays have the advantage of being  real-analytic, but
it doesn't appear that their origins
can be prescribed by this method. Moreover,
the relation of $R_j(T)$ to $F(T)$ is unclear.

\v

We now provide an outline of the paper.
The starting point  is the approximation
theorem for \K metrics by Bergman metrics:
For $k\geq 1$,
the space $\cH_k\sub\cH$  of Bergman metrics associated to $L^k$
is a finite dimensional symmetric Riemannian
sub-manifold.
If $h\in \cH$ and $h(k)\in \cH_k$ is the associated
Bergman metric, then the theorem of Tian-Yau-Zelditch \cite{Y},\cite{T97},\cite{Z}
implies $h(k)\ra h$ in the $C^\i$ topology.

\v
Now
fix $h_0, h_1\in\cH$, a pair of distinct elements, and let
$h(t;k)$ be the unique smooth geodesic segment in $\cH_k$  defined by the
conditions $h(0;k)=h_0(k)$ and  $h(1;k)=h_1(k)$.
It was proved in \cite{PS06}  that
the sequence $h(t;k)$ converges uniformly, in the weak $C^0$ sense of Theorem 1,
to a weak geodesic
segment
$h(t)$ in $\cH$ with the property: $h(0)=h_0$ and $h(1)=h_1$.  Moreover, $h(t)$  equals
the  $C^{1,1}$ geodesic joining $h_0$ to $h_1$, whose existence was established by
Chen~\cite{Ch}. We note that another approximation of the $C^{1,1}$ geodesic by
potentials $\tilde h(t;k)$ in $c_1(L)+{1\over k}c_1(K_X)$ has been very recently
constructed by Berndtsson \cite{Be06}.
\v

The proof of Theorem 1 follows the method of \cite{PS06}. First, we
construct a geodesic ray $h(t;k)=h_0e^{-\phi(t;k)}$ with $h(0;k)=h_0(k)$
that ``points
in the direction of $T$". Then we  prove that
\bea
\int_{X\times A}\O_k^{n+1}=O(k^{-1})
\eea
where $\O_k$ is associated to $\phi(t;k)$. This step relies  on the ideas
developed in
the recent work of Donaldson \cite{D05}. It also requires
some estimates on test configurations, which include the following very simple, but
basic estimate for the endomorphisms $A_k$ on $H^0(X_0,L_0^k)$ determined by
a test configuration,
\bea
\label{Ak}
\|A_k\|_{op} \ = \ O(k).
\eea
Next, we use the methods of pluripotential theory to establish the convergence of the
$\phi(t;k)$.
In the case of geodesic rays, the annulus $A$ is actually a punctured disk,
and the boundary behavior at the puncture has to be treated carefully, by
controlling the
asymptotics for the $\phi(t;k)$ at the puncture.

\v
For Theorem 2, we show that, when the test configuration is non-trivial,
the sup norm of $\phi_t$ goes to $\infty$ near the puncture.
This implies that the geodesic is non-trivial.
A key ingredient is Donaldson's formula \cite{D05} for the leading coefficient of
$\Tr(A_k^2)$.

\v
Theorem 3 is a direct consequence of the work of  Tian \cite{T97}
and Paul-Tian \cite{PT06}:
We apply the formula in  \cite{T97} which
relates the metric
of the CM line bundle $L_{CM}$ to the K-energy. We then use \cite{PT06}
which relates the line bundle $\lambda_{CM}$ on the Hilbert scheme to $L_{CM}$.

\v
We would like to add some references that have come
to our attention since the posting of the first version of
this paper.
In a paper \cite{PT06II} which appeared shortly after ours,
Paul and Tian present several results which include
in particular Theorem 3. In fact, they actually prove
a stronger result, in which the assumption on the injectivity
of the Chern map is removed.
As should be clear from its proof and as we already noted above,
Theorem 3 was in any case an immediate consequence of their earlier
work.
In the recent paper \cite{ChIII}, X.X. Chen shows
that geodesic rays parallel to a given geodesic ray can be
constructed under a certain assumption of tame ambient geometry.
 We would also like to note that constructions
involving upper envelopes appear frequently in
pluripotential theory, notably in the work of Kolodziej
\cite{K}.

\v\v

{\bf Acknowledgement:} We would like to thank Julius Ross for some
helpful conversations.

\newpage

\section{Test Configurations: preliminaries}
\setcounter{equation}{0}

\v
\subsection{Definition}

Let $L\ra X$ be an ample line bundle over a compact complex
manifold. A test configuration, as defined by Donaldson
 \cite{D02}, consists of the following
data:

\v

(1) A scheme $\cX$ with a $\C^\times $ action $\r$.

\v

(2) A $\C^\times$ equivariant line bundle $\cL\ra\cX$ which is ample on all fibers.

\v

(3) A flat $\C^\times $ equivariant map $\pi: \cX\ra \C$
where $\C^\times $ acts on $\C$ by multiplication
\v

satisfying the following: The fiber $X_1$ is
isomorphic to $X$ and the pair
$(X,L^r)$ is isomorphic to  $(X_1, L_1)$ where, for $w\in\C$, $X_w=\pi^{-1}(w)$
and $L_w=\cL|_{X_w}$. After raising ${\cal L}$ and $L$ to sufficiently high powers,
we may assume that $L$ is very ample, that $H^0(X,L)$ generates
$\oplus_{k=0}^\infty H^0(X,L^k)$, and that ${\cal L}$ has exponent one.
Thus we set $r=1$.

\v
If $\tau\in\C^\times$ and $w\in\C$, let
$\r_k(\tau,w):  H^0(X_w,L_w^k) \ra H^0(X_{\tau w},L^k_{\tau w})$
be the isomorphism induced by $\r$. If $w=0$ we write
$\r_k(\tau,0)=\r_k(\tau)$.
We also
let $B_k\in \End(V_k)$ be defined
by
\bea
\r_k(e^t)=e^{tB_k}
\eea
for $t\in\R$,
and $A_k $ the traceless part of $B_k$.
The eigenvalues of $A_k$ are denoted by $\l_0^{(k)}\leq\lambda_1^{(k)}\leq\cdots
\leq \lambda_{N_k}^{(k)}$, and the eigenvalues of $B_k$ are denoted by
$\eta_0^{(k)}\leq\eta_1^{(k)}\leq\cdots\leq\eta_{N_k}^{(k)}$.
Thus
$\r_k:\C^\times\ra GL(V_k)$
where $V_k=H^0(X_0,L_0^k)$.
Let $d_k=\dim V_k$ and $ w(k)=\Tr(B_k)$, the weight of the induced
action on $\det(V_k)$. Then, as was observed in
\cite{D02}, there is an asymptotic
expansion

\be
\label{expansion}
{ w(k)\over kd_k}\ = \ F_0 + F_1k^{-1}+F_2k^{-2} +\cdots \ \ \ \ {\rm as} \ \ k\to\i
\ee
The Donaldson-Futaki invariant $F(T)$, or simply Futaki invariant, of $T$ is defined by the formula:
$F(T)=F_1$.

\v

\subsection{Equivariant imbeddings of test configurations}

The construction of the Bergman geodesics associated to a
test configuration $T$ relies on the existence of an equivariant,
unitary imbedding of $T$ into projective space, whose
existence was first established by Donaldson \cite{D05}. In this
section, we begin by recalling the statement of Donaldson's result.

\v
Let $T$ be a test configuration of exponent $r=1$ for the pair $(X,L)$.
For $k$ large, since $L$ is very ample,
we have canonical compatible
imbeddings $\iota_k: X_1\sub \P(H^0(X_1,L_1^k)^*)$
and  $\iota_k:L^k_1\hookrightarrow O_1(1)$ where
$O_w(1)\ra \P(H^0(X_w,L_w^k)^*) $ is the hyperplane line bundle,
where $H^0(X_w,L_w^k)^*$ is the dual of $H^0(X_w,L_w^k)$.
\v
One can show that  the bundle $\pi_*\cL^k\ra\C$
has  an equivariant trivialization and thus the test
configuration has an equivariant imbedding into
projective space.
To be precise: Let $\T$ be an arbitrary
vector space isomorphism $\T: H^0(X_0,L_0^k)\ra H^0(X_1,L_1^k)$:
Let $\cX^\times = \pi^{-1}(\C^\times)$ and let $\cL^\times=\cL|_{\cX^\times}$. Define
an imbedding
$I_\T: (\cL^\times)^k\hookrightarrow O_0(1)\times \C^\times$ by the formula

\be I_\T(\r(\tau)l)\ = \ \left[(\r_k(\tau)\T^*(\iota_k(l) ), \tau\right]
\ee

where $\tau\in\C^\times$, $l\in L_1^k$ and $\T^*: O_1(1)\ra O_0(1)$ is
the isomorphism induced by the dual vector space
isomorphism $\T^*: H^0(X_1,L_1^k)^*\ra H^1(X_0,L_0^k)^*$. We similarly
define the imbedding $I_\T: \cX^\times \hookrightarrow \P(H^0(X_0,L_0^k))\times \C^\times$.
Then we say $\T:H^0(X_0,L_0^k)\ra H^0(X_1,L_1^k)$ is a ``regular generator of $T$"
if $I_\T$ extends to an imbedding
$\cL^k\hookrightarrow O_0(1)\times \C$ which restricts, over
the central fiber, to the
canonical embedding  $L_0^k\hookrightarrow O_0(1)$.

\v
Next let $h$ be a fixed metric
on $L$. It is shown in \cite{D05} that there exists an
regular generator $\T$
which respects $h$ structure
in the following sense: The metric $h$ defines a hermitian metric $H_k$ on $H^0(X,L^k)$
by the formula $\langle s,s'\rangle=\I_X (s,s')_{h^k}\ \o^n $ where
$\o$ is the curvature of~$h$. If $\T$ is a regular
generator of $T$, then we can use the isomorphism
$\T:V_k \ra H^0(X,L^k)$  to  define a metric on $V_k$, which
we  call $H_k(\T)$. Let $B_k$ be the endomophism of $V_k$
defined by: $\r_k(e^t)=e^{tB_k}$ for $t\in\R$.
We say $\T$ is a regular hermitian generator if $B_k$
is hermitian with respect to $H_k(\T)$. In other
words, $\T$ is regular hermitian if $\r_k(\tau):V_k\ra V_k$
is an isometry for $|\tau|=1$.

\v
In \cite{D05} the following is proved:
\begin{lemma}
\label{rigid}
Let $T$ be
a test configuration for $(X,L)$ and $h$
a positively curved metric on~$L$. Then
there exists $\T$, a regular hermitian generator
for $T$.
The metric $H_k=H_k(\T)$ is independent of the choice of
such a $\T$. Moreover,
the map $\T: V_k\ra H^0(X_1,L_1^k)$
is unique up to an isometry of  $V_k$ which commutes with $B_k$.
\end{lemma}

Our formulation of Lemma \ref{rigid} is somewhat
different than that given in  \cite{D05} and
in order to make the relationship between the two precise,
we
shall provide a complete proof (which
is of course essentially the one which appears
in \cite{D05}):

\v
Let $E\ra \C$ be an algebraic
vector bundle of rank $r$. Then
$E(\C)$, the space of global sections of $E$,
is a free $\C[t]$ module of rank $N+1$.
A ``trivialization of $E$'' is just a
choice of ordered basis $S_0,...,S_N$ of
the $\C[t]$ module $E(\C)$.
\v
If $S_0,...,S_N$ is a trivialization of $E$,
and if $t\in\C$, then $S_0(t),...,S_N(t)$
is a basis of the fiber $E_t$ so,
we have a well defined isomorphism
$\phi_{t_2,t_1}:E_{t_1} \approx E_{t_2}$
for any pair $t_1,t_2\in\C$,
which takes the basis $S_j(t_1)$ to
the basis $S_j(t_2)$. The collection $\{\phi_{t_2,t_1}\}$
defines a regular cocycle, that is:
$\phi_{t_3,t_2}\phi_{t_2,t_1}=\phi_{t_3,t_1}
$ and   for every $e\in E_{t_1}$,
the map $t\mapsto
\phi_{t,t_1}(e)$ is a global section of~$E$. Conversely,
a regular cocycle $\phi_{t_2,t_1}$ defines a trivialization of $E$.
\v
Now suppose $E\ra \C$ is a vector bundle with a
$\C^\times $ action, covering the usual action of
$\C^\times$ on $\C$. This means that we are given
an algebraic map $\r:\C^\times\ra \Aut(E\ra \C)$.  Thus,
if
$\tau\in\C^\times$ then $\r(\tau):E\ra E$ is a function
with the following properties:

\begin{enumerate}
\item The function $\r(\tau) $ maps the fiber $E_t$ into
the fiber $E_{\tau t}$, that is:
$\pi(\r(\tau)e)=
\r(\tau)\pi(e)
$.

\item The function $\r(\tau):E_t\ra E_{\tau t}$ is an
isomorphism of vector spaces.

\item If $\tau_1,\tau_2\in\C^\times$, then
$\r(\tau_1\tau_2)=\r(\tau_1)\r(\tau_2)$.

\item The map $\C^\times\times E\ra E $ given by
$(\tau,e)\ra \r(\tau)e $ is algebraic.
\end{enumerate}

\v
Let
$S_0,...,S_N$ be a basis of global sections for $E$. If $S:\C\ra E$
is an arbitrary global section, and if $\tau\in\C^\times$, then
$S^{\r(\tau)}(t)=\r(\tau)^{-1}S(\tau t)$ is also a
global section. Hence, there is a matrix $A(\tau,t)\in
GL(N+1,\C[\tau,\tau^{-1},t])$ with the property:
\be
\label{A}
 \underline{S}^{\r(\tau)}\ = \ A(\tau,t)\underline{S}
\ee
where $\underline{S}$ is the column vector whose components are
the $S_j$. Note that
$$ \underline{S}^{\r(\tau_2\tau_1)}\ = \ \r(\tau_2\tau_1)^{-1}
\underline{S}(\tau_2\tau_1 x)\ = \ \r(\tau_1)^{-1}
A(\tau_2,\tau_1 t)\underline{S}(\tau_1x)\ = \
A(\tau_2,\tau_1 t)A(\tau_1,t)\underline{S}(x)
$$
where, in the last equality, we are using the fact that
$\r(\tau_1)^{-1}$ is linear on the fibers. Hence:

\be
\label{cocycle}
 A(\tau_2\tau_1,t)\ = \ A(\tau_2,\tau_1 t)A(\tau_1,t)
\ee
In particular, if  $A(\tau)=A(\tau,0)$, then
$A(\tau):\C^\times\ra GL(N+1,\C)$ is a one parameter
subgroup.
\v
With these preliminaries in place, we now show that if $E\ra \C$ is
an vector bundle with $\C^\times$
action, then $E$ has a $\C^\times$
equivariant trivialization:

\begin{lemma}
\label{trivial equivariant}
Let $E\ra\C$ be a vector bundle of rank $r=N+1$ with a $\C^\times$ action. Then
there exists a basis of global sections $S_0,...,S_N$ such that
$A(\tau,t)$ is independent of $t$, that is,
$A(\tau,t)=A(\tau,0)\equiv A(\tau)$. In other words, there
exists a regular cocycle $\{\phi_{t_2,t_1}\}$ satisfying
\be
\label{equi phi}
\r(\tau)\phi_{t_2,t_1}\r(\tau)^{-1}\ = \ \phi_{\tau t_2,\tau t_1}
\ee
The basis  $S_0,...,S_N$ is unique
up to change of basis matrices $M(t)\in GL(N+1,\C[t]) $ with the
property:
$M(\tau t)= A(\tau)M(t)A(\tau)^{-1} $.
\end{lemma}

\v

{\it Proof.} Choose any $\C[t]$ basis $S_0,...,S_N\in E(\C)$  and define
$A(\tau,t)\in GL(N+1,\C[\tau,\tau^{-1},t]) $ as in
equation
(\ref{A}). Thus $\det(A(\tau,t))=a\tau^p =\det(A(\tau))$
for some integer $p$ and some $a\in\C^\times$. Now
consider the set

$$ \cS\ = \ \{ S_j^{\r(\tau)} : \tau\in\C^\times, 0\leq
j\leq N
\}
$$
 Let $V\sub E(\C) $ be the complex vector
space generated by $\cS$. We claim that $V$ is finite
dimensional and invariant under the action of
$\C^\times$. In fact, since $\underline{S}^{\r(\tau)}=A(\tau,t)
\underline S
$
we see that the $S_j^{\r(\tau)}$ are all linear
combinations, with $\C$ coefficients, of elements
in the set $\{ t^mS_j:
0\leq j\leq N, 0\leq m\leq M
\}
$, where $M$ is chosen so that the entries of
$A(\tau,t)$, which are polynomials in $t$ with
coefficients in $A[\tau,\tau^{-1}]$, all have degree
at most $M$.
\v
Choose a basis $\{T_\m; 0\leq \m\leq K\}$ of $V$ with the
property
$T_\m^{\r(\tau)}=\tau^{l_\m}T_\m $ for some integers
$l_\m$. Choose $\m_j, 0\leq j\leq N$, such that
$T_{\m_j}(0)$ are linearly independent. This can
certainly be done since the $T_\m$ span $V$, and $V$
contains the $S_j$. Let $\underline{T} $ be the column vector
consisting of the $T_{\m_j}$. Then $\underline T(t)=C(t)\underline S(t) $
for some $(N+1)\times (N+1) $ matrix $C(t)$ with coefficeints
in $\C[t]$, for which $C(0)$ is invertible.
The existence of such a matrix is guaranteed
by the fact that the $S_j$ form a $\C[t]$ basis of
$E(\C)$. Replacing $\underline S$ by $C(0)\underline S$ doesn't change
$V$ and allows us to assume $C(0)=I $. Now

$$ \underline T^{\r(\tau)}(t)\ = \ \r(\tau)^{-1} \underline T(\tau t)
\ = \ \r(\tau)^{-1}C(\tau t)\underline S(\tau t)\ = \ C(\tau
t)A(\tau,t)\underline S(t)
$$
On the other hand,
$\underline T^{\r(\tau)}(t)=U(\tau)\underline T (t)$ where $U(\tau)$ is
diagonal with diagonal entries of the form $\tau^l$.
Hence
$$U(\tau)\underline S(0)=U(\tau)\underline T(0)= T^{\r(\tau)}(0)=S^{\r(\tau)}(0)
=A(\tau)\underline S (0)
$$
 so $U(\tau)=A(\tau) $. Thus

$$
\ C(\tau t)A(\tau, t)\underline S(t)\ = \
\underline
T^{\r(\tau)}(t)\ =A\ (\tau)\underline T(t)\ = \ A(\tau)C(t)
\underline S(t)\
$$
which implies: $A(\tau)C(t)=C(\tau t)A(\tau, t) $. Since
$\det(A(\tau))=\det(A(\tau,t)) $ for all $t$, we have
$\det(C(\tau t))=\det(C(t))  $ which means that
$\det (C(t))$ is independent of $t$. Since $C(0)=I$, we
conclude $\det (C(t))=1$ and this implies that $\underline T$ is
a $\C[t]$ basis of $E(\C) $. This now establishes
Lemma \ref{trivial equivariant}.
\v
At this point we can prove the existence of
a regular generator for $T$:
Let $E=\pi_*(\cL^k)^*$ so that
$E_t=H(X_t,L_t^k)^*$. Then we
define $\T^*:E_1\ra E_0$ by the
formula: $\T^*=\phi_{0,1}$ where
$\phi_{t_2,t_1}$ satisfies
(\ref{equi phi}), with $\r(\tau)$ replaced
by $\r^*(\tau)=\r(\tau^{-1})^*$. One easily checks
that $\T$ is a regular generator of $T$.
\v

\begin{lemma}
\label{trivial hermitian} Let $H_1$ be a hermitian
metric on $E_1$. Then there is a unique equivariant
trivialization $\phi_{t_2,t_1}$ such that
$\r(\tau)^{-1}\phi_{\tau,1}:E_1\ra E_1 $ is an isometry
for  all $\tau\in\C^\times$ with $|\tau|=1$.

\end{lemma}

{\it Proof.} Let $\{\phi_{t_2,t_1}\}$ be
any equivariant trivialization.
Consider the decomposition $E_0=\oplus V_i $ into
eigenspaces for the action of $\C^\times$. Let $\tau^{w_j}$
be the restriction of $\r(\tau)$ to the subspace~$V_j$.
We may assume that $w_1<w_2<\cdots <w_l$. Thus
$\s_{j=1}^l w_j\dim(V_j)=N+1 =\dim(E_0)$. Let
$e_0,...,e_N $
of $E_0$ be given by the union of the
bases of the $V_j$ and define $S_j(t)=\phi_{t,0}(e_j)$
and let $W_i = \phi_{t,0}(V_i)\sub E_1$.
Then $S_0,...,S_N$ is a trivialization of $E\ra\C$.
\v
Let $A(\tau)$ be the diagonal matrix which represents
the automorphism
$\r(\tau): E_0\ra E_0$ with respect to the basis $e_j$.
Then $A(\tau)$ also represents the  automorphism
$\r(\tau)^{-1}\phi_{\tau,1}: E_1\ra E_1$ with respect
to the basis $S_j(1)$. We want to modify the
equivariant trivialization $\phi_{t_2,t_1}$ in
such a way that this automorphism is an isometry.
To do this, we must find a matrix
 $M(t)\in GL(N+1,\C[t])$ satisfying:

\begin{enumerate}

\item $M(\tau t)A(\tau)M(t)^{-1}=A(\tau) $ for all
$t,\tau$.

\item $M(1)S_j(1) $ is orthonormal with respect to $H$.

\end{enumerate}
The
first condition says that
$M(t)
$ is a block matrix with blocks $t^{w_i-w_j}\al_{ij} $
where
$\al_{ij}$ is independent of $t$. Since $M(t)\in GL(N+1,\C[t])$,
this implies that $\al_{ij}=0 $ if $i<j$. Thus
$M(t)$ is lower block triangular. On the other hand,
the usual Gram-Schmidt process allows us to choose
an $M(1)$ of this form which satisfies condition 2:
First choose an orthonormal basis of $W_0$. Then
choose an orthonormal basis of $W_0^\perp\sub W_0\oplus
W_1$, etc.
\v
Finally we prove uniqueness. Let $M(t)\in GL(N+1,\C[t])$
satisfy 1. and 2. and assume furthermore that the $e_j$
are orthonormal and that $M(0)=I$. Then we must
show that $M(t)=I$ for all $t$. Since the $e_j$ are
orthonormal, the matrix $M(1)$ is unitary. On the
other hand, it is lower block triangular. This
implies it is  block diagonal. Since the $i,j$
block is of the form $t^{w_i-w_j} \al_{ij }$, and
since $\al_{ij}=0$ for $i\not=j $, we see that $M(t) $
is independent of $t$ so $M(t)=M(0)=I $.
The lemma is proved.
\v

Note that if $\phi$ is any equivariant trivialization,
then $\r(\tau)^{-1}\phi_{\tau,1}:\C^\times \ra GL(E_1) $
is a homomorphism:

$$
\r(\tau_1)^{-1}\phi_{\tau_1,1}
\r(\tau_2)^{-1}\phi_{\tau_2,1}\ = \
\r(\tau_1\tau_2)^{-1}\phi_{\tau_1\tau_2,\tau_2}
\phi_{\tau_2,1}\ = \
\r(\tau_1\tau_2)^{-1}\phi_{\tau_1\tau_2,1}
$$
where the first equality makes use of the equivariance
property of $\phi$, and the second follows from
the cocycle property of $\phi$. Thus the theorem
can be restated as follows: There exists an equivariant
trivialization $\phi$ such that
$\r(\tau)^{-1}\phi_{\tau,1}: S^1\ra GL(E_1)$ is
a unitary representation.
\v
To deduce Lemma \ref{rigid} from Lemma \ref{trivial hermitian},
we again define $\T^*=\phi_{0,1}$ . Let $\tau\in\C^\times$ be of
unit length. Then to show $\r_k(\tau)^*:V_k^*\ra V_k^*$
is an isometry is equivalent, by definition of the metric on $V_k$, to
showing $(\T^*)^{-1}\r_k(\tau)^*\T^*: H(X_1,L^k)^*\ra H(X_1,L^k)^*$
is an isometry. Thus we must show
$\phi_{1,0}\r_k(\tau)^*\phi_{0,1}=\phi_{1,0}\r_k^*(\tau^{-1})\phi_{0,1}$
is an isometry.
But (\ref{equi phi}) implies
$\phi_{1,0}\r_k^*(\tau^{-1})\phi_{0,1}=
\r^*(\tau)^{-1}\phi_{\tau,0}\phi_{0,1}=\r^*(\tau)^{-1}\phi_{\tau,1}$
which is an isometry by the result of Lemma  \ref{trivial hermitian}.
This proves Lemma 1.
\v



\section{Estimates for test configurations}
\setcounter{equation}{0}

\subsection{Bounds for $A_k$}

Let $T$ be a test configuration and define the endomorphisms $A_k$ and $B_k$
and their eigenvalues $\lambda_\al^{(k)}$ and $\eta_\al^{(k)}$ as in Section \S 2.1.
The following simple estimate for the operator norm $\|A_k\|_{op}$ of the endomorphisms
$A_k$ plays an important role in the subsequent bounds for the total masses
of the Monge-Amp\`ere currents:

\v
\begin{lemma}
\label{bound}
There is a constant $C>0$ which is independent of $k$ such that
$|\l^{(k)}_\al|\leq Ck$ for all $k>0$ and all $\al$ such that
$0\leq \al\leq N_k$.
\end{lemma}

{\it Proof of Lemma \ref{bound}}.
After applying Lemma \ref{rigid} with $k=1$,
we may assume that $\cX\sub\P^m\times \C$, $m=N_1+1$, and that
$\r(\tau)$ is a diagonal matrix in $GL(m+1)$ whose entries are
$\tau^{\eta_0},...,\tau^{\eta_m}$ where
$\eta_0\leq \cdots \leq \eta_m$ are integers. The scheme $X_0\sub\P^m$ is
defined by a homogenous ideal $I\sub \C[X_0,...,X_m]$ and
we write

$$ \C[X_0,...,X_m]/I\ = \ \bigoplus_{k\geq 0} S_k/I_k
$$
where $S_k\sub \C[X_0,...,X_m]$ is the space of polnomials
which are homogeneous of degree $k$ and $I_k=S_k\cap I$. Then,
for $k>>0$, we have $H^0(X_0,L_0^k)= S_k/I_k$. The matrix
$\r(\tau)$ defines an automorphism of $\C[X_0,...,X_m]$,
determined by the formula: $X_j\mapsto \tau^{\eta_j}X_j$.
This automorphism
leaves $S_k$ and $I_k$ invariant, and thus it induces an
automorphism of $S_k/I_k$ which is, by definition, the
map $\r_k(\tau)$.
\v
The monomials of degree $k$ form a basis of $S_k$ which
are eigenfuntions of $\r(\tau)$: More precisely, if $X^p$
is a monomial, with $p=(p_0,...,p_m)$ and $p_0+\cdots +p_m=k$,
then we have $\r(\tau)\cdot X^p= \tau^{p\cdot \eta}X^\al$.
Since the monomials of degree $k$ span $S_k/I_k$, some subset
form a basis of eigenvectors for that space. Thus the eigenvalues
of the $B_k$ form a subset of $\{p\cdot\eta: p_0+\cdots+p_m=k\}$.
On the other hand, for such an $p$, we clearly have
$|p\cdot \eta|\ \leq \sup|\eta_j|\cdot k $ and this proves
that
\bea
|\eta_\al^{(k)}|\ \leq \ C\,k,
\eea
with $C=\sup_{0\leq j\leq m}|\eta_j| $. On the other hand,
\bea
\lambda_\al^{(k)}=\eta_\al^{(k)}-{{\rm Tr}\,(B_k)\over N_k+1}
=\eta_\al^{(k)}+O(k).
\eea
This proves Lemma \ref{bound}.

\v

\v

\subsection{An alternative characterization of the Futaki invariant}

\subsubsection{The $F_{\o}^0$ functional}

Let $X$ be a compact complex manifold of dimension $n$ and $\o=\o_0$ a \K metric on $X$. Let
$\cH=\cH_{\o}$ be the space of \K potentials:
\be
\cH_\o\ = \ \{\phi\in C^\i(X): \o_\phi=\o+ {\sqrt{-1}\over 2}\ddb\phi>0\}
\ee
The functionals $F^0_\o,\,\nu_\o: \cH\ra \R$ play an important role in
\K geometry and are defined as follows:

\bea
F^0_\o(\phi)  \ &=& \ \ -{1\over n+1}\left(\I_X\o^n\right)^{-1} \s_{j=0}^{n}\ \I_X
\phi\o_\phi^j\o^{n-j} \ = \ -{1\over n+1}\left(\I_X\o^n\right)^{-1}E_\o(\phi)
\nonumber\\
\nu_\o(\phi) \ &=&-\left(\int_X\o^n\right)^{-1}\int_0^1\int_X\dot\phi(s-\hat s)\o_t^n\,dt,
\eea
Here $\phi_t$, $0\leq t\leq 1$, is a smooth path in ${\cal H}_\o$ joining the potential $\phi_0$
for $\o_0$ to $\phi=\phi_1$. Then a simple calculation shows

\be
\label{derivative}
\dot E_\o(\phi_t)\ = \ (n+1)\I_X\dot\phi_t\ \o_{\phi_t}^n\ \ \ {\rm and} \ \ \
\ddot E_\o(\phi_t)\ = \ (n+1)\I_X (\ddot\phi_t - |\pl\dot\phi_t|^2)\ \o_{\phi_t}^n
\ee
Thus $E$ satisfies the cocycle property: $E_\o(\phi)+E_{\o_\phi}(\psi)=E_\o(\phi+\psi)$\ .
Note as well that if $f:Y\ra X$ is a biholomorphic map, then
\be\label{biholomorphic}
E_{f^*\o}(\phi\circ f)\ = \ E_\o(\phi)
\ee

\subsubsection{The Chow weight and the Futaki invariant}

Let $V$ be a finite dimensional vector space,
$Z\sub \P(V)$ a smooth subvariety, and $B\in gl(V)$.
Then we wish to define
the generalized Chow weight
$\m(Z,B)\in\R$. We start by assuming the $V=\C^{N+1}$ so that
$B$ is a  $(N+1)\times (N+1)$ matrix. Let
$\o_{FS}$ be the Fubini-Study metric on $\P^N$. We shall also denote by
$\o_{FS}$ the restriction of the Fubini-Study metric to $Z$. For $t\in\R$
let $\sigma_t\in GL(N+1,\C)$ be the matrix $\sigma_t=e^{tB}$ and
let $\psi_t:\P^N\ra\R$ be the function

\be
\label{psi}
\psi_t(z)\ = \ \log{|\sigma_t z|^2\over |z|^2}
\ee

Here we view $z$ as an element in $\P^N$ and, when there is no
fear of confusion, a column vector in $\C^{N+1}$.
\v
Then $\psi_t$ is a smooth path in $\cH$: In fact,
$\sigma_t^*\o_{FS}=\o_{FS}+{\sqrt{-1}\over 2}\ddb\psi_t$.
Define

\be
\label{mu}
\m(Z,B)\ = \ -\lim_{t\ra -\i}\dot E_{\o_{FS}}(\psi_t)\ = \ - \dot E(-\i)
\ee

Note that the
function $E(t)=E_{\o_{FS}}(\psi_t):\R\ra\R$
is  convex (see \cite{PS02, PS02a}), so the limit in (\ref{mu})
exists.

\v

Next we compute the derivative of $E(t)$:

\be
\label{dot}
{d\over dt}E_{\o_{FS}}(\psi_{t})\ =
\ {(n+1)}\I_Z {
{z^*\si_t^*}\cdot (B+B^*)\cdot \si_tz \over
z^*\si_t^*\si_tz}\ \sigma_t^*\o_{FS}^n\ = \
{(n+1)}\I_{\si_t(Z)} {
{z^*}\cdot (B+B^*)\cdot z \over
z^*z}\ \o_{FS}^n
\ee
where, for  $C$  a matrix with complex entries, we write $C^*={^t}\bar C$.
In particular,
\be
\label{dot zero}
 \dot E_{\o_{FS}}(\psi_t)|_{t=0}\ = \ \dot E(0)\ = \ {(n+1)} {\rm Tr}((B+B^*)\cdot M)
\ee
where
\be
\label{M} M_{\al\b}\ = \ M_{\al\b}(Z)\ = \ \I_Z {z_\al\bar z_\b\over \|z\|^2}\ \o_{FS}^n
\ee
\v

\begin{lemma}
\label{welldefined}
Let $V$ be a finite dimensional complex vector space, $B\in gl(V)$ and
$Z\sub \P(V)$ a smooth subvariety. Let $\t: V\ra\C^{N+1}$
be an isomorphism. Then $\m(\t(Z), \t B \t^{-1})$ is
independent of~$\t$.
\end{lemma}

{\it Proof.}
 We make use of the formula of Zhang \cite{Z96} and  Paul \cite{P00}
(see also \cite{PS02}):
If $Z\sub\P^N(\C)$ is a subvariety of dimension $n$ and degree $d$,
let ${\rm Chow}(Z)\in \P(H^0(Gr(N-n,\C^{N+1}, O(d))))$ be the
Chow point of $Z\sub\P^N$.
If $B\in gl(N+1,\C)$, $\sigma_t=e^{tB}$, and
$\psi_{\si_t}=\log{|\sigma_t(z)|^2\over |z|^2}$ then

\be
\label{zhang}
E_{\o_{FS}|_Z}(\psi_{\sigma_t})\ = \ \log{\|\sigma_t\cdot {\rm Chow}(Z)\|^2
\over \|{\rm Chow}(Z)\|^2}
\ = \
\ \log{\| {\rm Chow}(\sigma_tZ)\|^2\over \|{\rm Chow}(Z)\|^2}
\ee
\v
where $\|\cdot \|$ is the Chow norm defined on $H^0(Gr(N-n,\C^{N+1}, O(d)))$.
\v
Suppose $M\in GL(N+1,\C)$. Then
\be
\label{zhang1}
E_{\o_{FS}|_{MZ}}(\psi_{M\sigma_t M^{-1}})\ = \ \log{\|M\sigma_tM^{-1}\cdot
{\rm Chow}(MZ)\|^2\over \|{\rm Chow}(MZ)\|^2}
\ee
Subtracting (\ref{zhang}) from (\ref{zhang1}) we get

$$ E_{\o_{FS}|_{MZ}}(\psi_{M\sigma_t M^{-1}}) \ - \ E_{\o_{FS}|_Z}(\psi_{\sigma_t})
\ = \ \log{\|M\sigma_t\cdot {\rm Chow}(Z)\|^2\over \|\sigma_t\cdot{\rm Chow}(Z)\|^2}
\ -\ \log{\|M\cdot {\rm Chow}(Z)\|^2\over \|{\rm Chow}(Z)\|^2}
$$
which is a bounded function of $t$, and hence the limit of its first
derivative is zero. This proves Lemma \ref{welldefined}.

\v
Now let $Z\sub \P(V)$ and $B\in gl(V)$.
Let $\t:V\ra\C^{N+1}$ be an isomorphism and define $\m(Z,B)=\m(\t(Z), \t B
\t^{-1})$. The lemma guarantees that this definition is unambiguous.
Note that (\ref{zhang}) shows that $\m(Z,B)$ is just the usual Chow
weight. (The Chow weight is normally defined only  when $B$ is
a traceless diagonalizable matrix with integer eigenvalues,
but we find it convenient to work with this somewhat more
general notion).
\v

If $\tau\in GL(V)$ then
$\m(\tau(Z), B)= \m(\t\tau(Z), \t B \t^{-1})=\m((\t\tau)(Z),
(\t\tau)\tau^{-1}B\tau(\t\tau)^{-1})$.
We conclude
\be
\label{cocycle} \m(\tau(Z),B)\ = \ \m(Z, \tau^{-1}B\tau)
\ee
In particular, if $\tau$ commutes with $B$, then $\m(Z,B)=\m(\tau(Z),B)$.
\v
If we replace the functional
$E$ by $\n$, the $K$ energy functional,
we may define a corresponding invariant
$\tilde\m(Z,B)$ for $Z\sub\P^N(\C)$ and
$B\in gl(N+1,\C)$:
\be
\tilde \m(Z,B)\ = \ \lim_{t\to-\i}
\dot\n_{\o_{FS}}(\psi_t)
\ee
\v

It will be convenient for us to introduce an alternative characterization
of the Futaki invariant: Fix, once and for all,
an isomorphism $\kappa: (X,L^r)\ra (X_1, L_1)$.
We continue to assume that $r=1$ (the case $r>1$ can be
treated in a similar fashion).
Then  we have an induced isomorphism \ $H^0(X, L^k) = H^0(X_1,L_1^k)$.

\v

Let $\T$ be an equivariant trivialization
of $\pi_*\cL^k$. Then
 $I_\T|_{X_1}: X_1\hookrightarrow \P(H^0(X_0,L_0^k))$.
Let $Z_k\sub\P(H^0(X_0,L_0^k)^*)$ be the image of $I_\T|_{X_1}$
and $Z_k^{(0)} $ the image of the canonical imbedding
$X_0\sub \P(H^0(X_0,L_0^k)^*)$.
Note that $Z_k$ depends on the choice of $\T$, but that if $\T'$ is
another choice, then $\T'= U\T$ where $UA_k=A_kU$,
and thus the value $\m(Z_k,A_k)$ is independent of the choice of
equivariant $\T$

\begin{lemma}
\label{fut}
We have
\be
F(T)\ =  \ -c(X,\o)\cdot\lim_{k\to\i}{\m(Z_k,A_k)\over k^n}
\ee
where $c(X,\o)= {1\over n!(n+1)!}\I_X\o^n$.
\end{lemma}

{\it Proof. } Since this argument is implicit in Donaldson \cite{D02},
we only briefly sketch the proof (see as well Ross-Thomas \cite{RT}):
If $Z\sub\P^N$ and $\l:\C^\times\ra SL(N+1,\C)$ is a
one parameter subgroup, let $A\in sl(N+1)$ be such
that $\l(e^t)=e^{tA}$ and $Z^{(0)}=\lim_{\tau\ra 0}\l(\tau)(Z)$ (the flat
limit) so  $Z^{(0)}\sub \P^N$ is a subscheme of $\P^N$ with
the same Hilbert polynomial as~$Z$. Let $M_0=O(1)|_{Z^{(0)}}$.
Then $\l(\tau)$ defines an automorphism of $H^0(Z^{(0)},M_0^p)$
and we let $\tilde w(Z, A, p)$ be the weight of this action
on $\det(H^0(X_0,M_0^p))$. It is known that $\tilde w(p)$ is a polynomial
in $p$ for $p$ large such that

\be\label{tildew} \tilde w(Z,A, p) \ = \ {\m(Z,A)\over (n+1)!}\cdot p^{n+1}\ + \ O(p^n)
\ \ {\rm and} \ \ \tilde w(Z^{(0)},1)=0
\ee
(see for example, Mumford \cite{Mumford}).
Now let $T$ be a test configuration, let $r>0$ and consider
$Z_r\sub \P(H^0(X_0,L_0^r)^*)$.
Applying (\ref{tildew}) to $Z=Z_r$, $A=rN_rA_r$ and $M_0=L_0^r$, we get

\be\label{tildew1}
 \tilde w(Z_r, rN_rA_r, p) \ = \ {\m(Z_r,rN_rA_r)\over (n+1)!}\cdot p^{n+1}\ + \ O(p^n)
\ee
On the other hand, since $M_0^p=L_0^{rp}$, we get, with $k=rp$:

\be \tilde w(Z_r, rN_rA_r, p) \ = w(k)rN_r-w(r)kN_k\ = \ e_T(r)k^{n+1}\ + \ O(k^n)
\ee
where $e_T$ is a polynomial in $r$ of degree at most $n$.
If follows from the definition of $F(T)$ that
$-F(T)$ is the leading coefficient of $e_T(r)$. Comparing with
(\ref{tildew1}) we get

$$\lim_{r\to\i}\ {\m(Z_r,rN_rA_r)\over r^nr^{n+1} (n+1)!}\ = \ -F(T)
$$
Since $r^{-n}N_r = {1\over n!}\I\o^n + O(r^{-1})$, Lemma \ref{fut} follows.

\v

\section{Completion of the Proof of Theorem 1}
\setcounter{equation}{0}

\subsection{The Tian-Yau-Zelditch expansion}
Let $L\ra X$ be an ample line bundle over a compact
complex manifold $X$. If $h$ is a smooth hermitian
metric on $L$ then the curvature of $h$ is given
by $\o=R(h)=-{\sqrt{-1}\over 2}\ddb\log h$.
Let $\cH$ be the space of positively
curved hermitian metrics on $L$.
Then  $\cH$ contains a canonical family of
finite-dimensional negatively curved symmetric
spaces
$\cH_k$, the space of Bergman metrics, which  are defined as follows: For $k>>0$
and for
$\us=(s_0,....,s_{N_k})$ an ordered basis of
$H^0(X,L^k)$, let
$$ \iota_\us: X\hookrightarrow \P^{N_k}
$$
be the Kodaira imbedding given by
$x\mapsto (s_0(x),..., s_{N_k}(x))$. Then we have a
canonical isomorphism
$\iota_\us: L^k \ra \iota_\us^*O(1)$ given by

\be
\label{canonical}
\iota_\us(l)\ = \  \Big[({s_0\over s},{s_1\over s}, ... , {s_N\over s})
\mapsto {l\over s} \Big]
\ee
where $l\in L^k$ and $s$ is any locally trivializing section of $L^k$.
\v
Fix $h_0\in\cH$. Let $h_{FS}$
be the Fubini-Study metric on $O(1)\ra \P^{N_k}$ and
let
\be\label{fubini}
h_\us=
(\iota_\us^*h_{FS})^{1/k}\
=
\  {h_0\over
\left(\s_{\al=0}^{N_k} |s_\al|^2_{h_0^k}\right)^{1/k}}
\ee
Note that the right side of (\ref{fubini}) is
independent of the choice of $h_0\in \cH$. In
particular
\be\label{fubini two}
\s_{\al=0}^{N_k} |s_\al|^2_{h_\us^k}\ = \ 1
\ee
Let
$$ \cH_k\ = \ \{h_\us: \us \hbox{\ a basis
of  $H^0(X,L^k) $\} }\ \sub \ \cH
$$

\v
Then
$\cH_k = GL(N_k+1)/U(N_k+1)$ is a finite-dimensional
negatively curved symmetric space
sitting inside of $\cH$.  It
is well known that the $\cH_k$ are  topologically
dense in $\cH$: If $h\in \cH$ then there
exists
$h(k)\in\cH_k$ such that $h(k)\ra h$ in the
$C^\i$ topology. This follows from the Tian-Yau-Zelditch
theorem  on
the
density of states (Yau \cite{Y87}, Tian \cite{T90} and
Zelditch \cite{Z}; see also Catlin
\cite{C} for corresponding
results for the Bergman kernel). In fact, if
$h\in\cH$, then
there is a canonical choice of the approximating
sequence $h(k)$:
Let $\us$ be a basis
of $H^0(X,L^k)$ which is orthonormal with respect
to the metrics $h$.  In
other words,
\be
\label{metric}
\langle s_\al,s_\b\rangle_h\ = \ \I_X (s_\al,s_\b)_{h^k} \
\o^n=\d_{\al\b}\ \ \hbox{where $\o=R(h) $} \ .
\ee
The basis $\us$ is unique up to an element of $U(N_k+1)$.
Define $\r_k(h)=\r_k(\o) =\s_\al |s_\al|^2_{h^k}$.
Then  Theorem 1 of \cite{Z}, which is the $C^\i$ version
of the $C^2$ approximation result first estsablished
in  \cite{T90}, says
that for
$h$ fixed, we have a $C^\i$ asymptotic expansion as
$k\to\i$:
\be\label{zelditch} \r_k(\o)\ \sim \ k^n +
A_1(\o)k^{n-1} + A_2(\o)k^{n-1}+\cdots
\ee
Here the $A_j(\o)$ are smooth functions on $X$
defined locally by $\o$ which can be computed in terms
of the curvature of $\o$ by the work of Lu \cite{Lu}. In particular,
it is shown there that
\be A_1(\o)={s(\o)\over 2\pi}
\ee
 where $s(\o)$ is the scalar
curvature of $\o$.
\v
 Let
$\hat\us=k^{-n/2}\us$ and $h(k)=h_{\hat\us}$. Then
(\ref{fubini}) and (\ref{zelditch}) imply that
\be\label{powers of k}
{h(k)\over
h}\ = \ 1-{s(\o)\over 2\pi}\cdot{1\over k^2}+O\left({1\over k^3}\right)\ \ , \ \
\o(k)=\o + O\left({1\over
k^2}\right)\ \ , \ \
\phi(k)=\phi + O\left({1\over
k^2}\right)
\ee
Here, as before, $\o=R(h)$,
$\o(k)=R(h(k))$, $h=h_0e^{-\phi} $ and
$h(k)=h_0e^{-\phi(k)}$. In particular, $\o_0+{\sqrt{-1}\over 2}\ddb\phi(k)=
\o(k)={1\over k}\iota_\us^*\o_{FS}$.

\v

Lemma \ref{rigid} can now  be conveniently reformulated as follows:
\v
\begin{lemma}
\label{reformulation}
Let $\r:\C^\times \ra \Aut(\cL\ra\cX\ra \C)$ be a
test configuration $T$ of exponent one for the pair $(X,L)$, where
$L\ra X$ is  ample. Let $h_0$ be a positively curved metric
on $L\ra X$. Let $k$ be an integer such that $L^k$ is very ample. Then there
is
\v
(1) An orthonormal basis $\us=(s_0,...,s_{N_k})$ of $H^0(X,L^k)=H^0(X_1,L_1^k)$,

(2)
An imbedding
$I_{\us}: (\cL^k\ra\cX\ra \C)\hookrightarrow (O(1)\times \C\ra \P^{N_k}\times \C \ra \C)$

\v
satisfying
the following property: the imbedding $I_\us$ restricts to $\iota_\us$ on the fiber $L_1^k$
and $I_\us$ intertwines $\r(\tau)$ and $\tau^{B_k}$. More precisely: for every $\tau\in
\C^\times$ and every
$l_w\in
L_w^k$,

\be\label{I} I_\us(\r(\tau)l_w)\ = \ \left(\tau^{B_k}\cdot I_\us\left(l_w\right), \tau w\right)
\ee
where $\tau^{B_k}$ is a diagonal matrix whose eigenvalues are the
eigenvalues of $\r_k(\tau):V_k\ra V_k$.
\v
The matrix $B_k$ is uniquely determined, up to a permutation of the diagnonal
entries, by $k$ and the test configuration $T$.
Moreover, the basis $\us$ is uniquely determined by $h_0$ and $T$, up to an element of $U(N_k+1)$
which commutes with $B_k$. The image of $X_1$ is $Z_k\sub\P^{N_k}$.

\end{lemma}

\v

\subsection{Growth bounds for the Bergman geodesic rays}

We make precise the notation which appears in Theorem 1:
Let $L\ra X$ be an ample line bundle over a compact complex
manifold, and $\cH$ the space of positively curved metrics
on $L$. Let $h_0\in\cH$ and let $T$ be a test configuration for the pair $(X,L)$
of exponent $r$. We wish to associate to the pair $(h_0,T)$
an infinite geodesic ray in $\cH$ whose initial point is $h_0$.
After replacing $L$ by $L^r$ we may assume, without loss of
generality, that $r=1$ and that $L$ is very ample.

\v
Let $k$ be a large positive integer and choose
 $\us$, an orthonormal basis of $H^0(X,L^k)$ as
in Lemma \ref{reformulation}. Define $A_k$ to be the traceless part of $B_k$
and let $\l^{(k)}_0\leq\l^{(k)}_1\leq \cdots \leq \l^{(k)}_{N_k}$ be the diagonal
entries of $A_k$. Set $\hat \us=k^{-n/2}\us$
so that $h_{\hat\us}= h_0(k)$, where $h_{\hat\us}$ is
defined as in (\ref{fubini}).
Now let $\hat\us(t;k)= (e^{t\l_0}\hat s_0, e^{t\l_1}\hat s_1, ..., e^{t\l_{N}}\hat s_N)$,
and define

\be
h(t;k) \ = \ h_{\hat\us(t;k)}\ = \ h_0e^{-\phi(t;k)}\ = \ h_0(k)e^{-(\phi(t;k)-\phi(k))}
\ee
so that $h(t;k):(-\i,0]\ra \cH_k$ is a geodesic ray in $\cH_k$ and $h(0;k)=h_0(k)$.
In particular we have

\be\label{phit}
\phi(t;k)\ =  {1\over k}\log\left(k^{-n}\cdot\sum_{\al=0}^{N_k}\ e^{2t\l_\al}|s_\al|^2_{h_0^k}\right)
\ = \
{1\over k}\log\left(k^{-n}\cdot\sum_{\al=0}^{N_k}\ e^{2t\l_\al}|s_\al|^2_{h_0(k)^k}\right)\ + \ \phi(k)
\ee

\v
Let
\be
f(k)= {w(k)\over kd_k}-F_0\ = \ {F(T)\over k}+ O({1\over k^2})
\ee
where $w(k), d_k$ and $F_0$ are defined as in (\ref{expansion}).
In particular, $f(k)=O({1\over k})$.
\begin{lemma}
\label{phi bound lemma}
Let $k,l$ be positive integers with $k<l$. Then there exists $C_{k,l}>0$ with the
following property:

\be\label{phi bound}
-C_{k,l} \ < \ [\phi(t;l)+2t\cdot f(l)]\ - \
[\phi(t;k)+2t\cdot f(k)]\
< \ C_{k,l}
\ee

\end{lemma}

\v
{\it Proof.} If suffices to prove (\ref{phi bound}) in the case $k=1$.
Then, replacing $l$ by $k$, we have

\be\label{tilde phi}
\tilde \phi(t;k)-\tilde\phi(t;1)=
[\phi(t;k)+2t\cdot f(k)]\ - \
[\phi(t;1)+2t\cdot f(1)]\ = \
\log{
\left(k^{-n}\cdot\sum_{\al=0}^{N_k}
\ e^{2t\eta_\al^{(k)}}|s^{(k)}_\al|^2_{h_0^k}\right)^{1/k}
\over
\left(\sum_{\b=0}^{N}
\ e^{2t\eta_\b^{(1)}}|s_\b|^2_{h_0}\right)
\\ \\ }
\ee
where $\eta^{(k)}_0\leq \eta^{(k)}_1\cdots \leq \eta^{(k)}_{N_k}$
are the eigenvalues of the diagonal matrix $B_k$, $N=N_1$
and $s_\b=s_\b^{(1)}$.
\v
We now have

\be
\log(\sum_{\b=0}^{N}
\ e^{2t\eta_\b^{(1)}}|s_\b|^2_{h_0})\ = \
{1\over k}\log(\sum_{\al=0}^{N_k}
 e^{2t\eta_\al^{(k)}}|\tilde s^{(k)}_\al|^2_{h_0^k})\ + \ O(1)
\ee
where the $O(1)$ term is independent of $t$, and, for $\eta\in {\bf Z}$,

\be
\{\tilde s_\al^{(k)}:\eta_\al^{(k)}=\eta\}
\subset\{ s_0^{p_0}\otimes \cdots \otimes s_{N}^{p_N}\in H^0(X,L^k): \s_\b p_\b=k\ \
{\rm and} \ \
\s_\b p_\b\eta_\b^{(1)}\ = \ \eta\},
\ee
is a maximally linearly independent subset. On the other hand,
$(s_0^{(k)},...,s_{N_k}^{(k)})$ and
$(\tilde s_0^{(k)},...,\tilde s_{N_k}^{(k)})$ are two
bases of the same vector space which differ by
a lower block triangular matrix.
This proves Lemma \ref{phi bound lemma}.
\v

\subsection{The volume formula}

Let $\phi_t: [a,b]\ra\cH_\o$ be a smooth
path and let  $U_{a,b}=\{w\in\C^\times: e^a\leq |w|\leq e^b\}$. Let
$M_{a,b}=X\times U_{a,b}$ and $\O_{0}$ be the $(1,1)$ form on $M_{a,b}$
defined by pulling back $\o_0$. Define $\Phi(z,w): M_{a,b} \ra \R$  by

\be
\label{11form}
\Phi(z,w)\ = \ \phi_t(z)\ \ \ {\rm where} \ \ \ t=\log |w|
\ee

Let $\O_\Phi$ be the $(1,1)$ form on $M_{a,b}$ defined by
$\O_\Phi =\O_0+{\sqrt{-1}\over 2}\ddb\Phi$.
Then
\be
\O_\Phi^{n+1}\ = \ {1\over 4}(\ddot\phi_t\ - \ |\pl\dot\phi_t|^2)\o_{\phi_t}^n
\wedge ({\sqrt{-1}\over 2}dw\wedge d\bar w)
\ee
In particular, we have the key observation of  \cite{M87}, \cite{S92}, \cite{D99} :
\be \O_\Phi^{n+1} = 0 \ \iff \ \ddot\phi_t\ - \ |\pl\dot\phi_t|_{\o_{\phi_t}}^2=0
\ \iff \ \phi_t \ \hbox{is a smooth geodesic in $\cH$}
\ee
We say that a function $\phi_t(x)$ on $[a,b]\times X$ is a weak geodesic if
$\Phi$ is bounded, plurisubharmonic with respect to $\O_0$, and if $\O_\Phi^{n+1}=0$.
\v
Finally,  we obtain, using (\ref{derivative}), the following useful volume formula
\cite{PS06}:
\be
\label{volume}
{( n+1)}\I_{X\times U_{a,b}} \O_\Phi^{n+1}\ = \ \dot E_\o(b)\ - \ \dot E_\o(a)
\ee
where $E_\o(t)=E_\o(\phi_t)$.

\subsection{Volume estimates for the Monge-Amp\`ere measure}

We first need a few lemmas: Let $D^\times=\{w\in\C: 0<|w|<1\}$. We associate to $\phi(t;k)$
the function  $\Phi(k)$ on $X\times D^\times$ as in (\ref{11form}):

\be
\Phi(k)(z,w)\ = \ \phi(t;k)(z)\ \ \ {\rm where} \ \ \ t=\log |w|
\ee
and we let $\O_0$ be the pullback of $\o_0$ to $X\times D^\times$ and we let
$\O_{\Phi(k)}=\O_0+{\sqrt{-1}\over 2}\ddb\Phi(k)$.
\begin{lemma}
\label{volest}
We have $\lim_{k\to\i}\I_{X\times D^\times}\O^{n+1}_{\Phi(k)}=0$. In fact,
\end{lemma}

\be
\label{vol2}\I_{X\times D^\times } \O_{\Phi(k)}^{n+1}\ = \ O({k^{-1}})
\ee
{\it Proof.}
According to (\ref{volume}),
\be
\label{v1}\I_{X\times D^\times } \O_{\Phi(k)}^{n+1}\ =
 \ \I_X\dot \phi(0;k)\ \o_{\phi(0;k)}^n\ - \
\lim_{t\ra\-\i}\I_X \dot\phi(t;k)\ \o_{\phi(t;k)}^n
\ee
Hence it suffices to show that each of the two terms in (\ref{vol2}) is $O({1\over k})$.

\v

Let $\psi(t;k)= k[\phi(t;k)-\phi(k)] + n\log k$. Then  $\psi(t;k)=\log{|\sigma_t z|^2\over |z|^2} $
where $\sigma_t=e^{tA_k}$. Recall that $\o_0+ {\sqrt{-1}\over 2}\ddb\phi(k)=\o_0(k)$
and $\o_0(k)={1\over k}\iota_\us^*\o_{FS}$.
Thus
$$\o_{\phi(t;k)}\ =\
\o_0+{\sqrt{-1}\over 2}\ddb\phi(t;k)\ =\
\o_0(k)+{1\over k}{\sqrt{-1}\over 2}\ddb \psi_t
\ =\ {1\over k}\iota_\us^*\left(\o_{FS}+{\sqrt{-1}\over 2}\ddb\psi_t\right)
$$
Since we also have $\dot\phi(t;k)= {1\over k}\dot\psi(t;k)$ we conclude:
\be (n+1)\I_X \dot\phi(t;k)\ \o_{\phi(t;k)}^n\
\ = \
(n+1)\cdot{1\over k}\cdot{1\over k^n}\cdot\I_{Z_k} \dot\psi_t \
\left(\o_{FS}+{\sqrt{-1}\over 2}\ddb\psi_t\right)^n \
=
\ {1\over k}\cdot{1\over k^n}\cdot \dot E_{\o_{FS}}(\psi_t)
\ee

Thus

\be -\lim_{t\to-\i}\I_X \dot\phi(t;k)\ \o_{\phi(t;k)}^n\ =
\ {1\over (n+1)}{1\over k}{\m(Z_k,A_k)\over k^n}
\ee

Now, according to Lemma \ref{fut}, ${\m(Z_k,A_k)\over k^n}$
has a finite limit as $k$ tends to infinity. In fact, the limit
is equal to $F(T)$, the Futaki invariant. Thus
the second term in (\ref{v1}) is $O({1\over k})$.

\v

In order to treat the first term, we require the following
result from \cite{D05}:

\v

\begin{lemma}
\label{D05}
Let $L\ra X$ be an ample line bundle over a compact complex
manifold $X$ and $h$ a metric on $L$ with positive curvature $\o$.
Let $\us$ be an orthonormal basis of $H^0(X,L^k)$ and
let $\iota_{\hat\us}: X\hookrightarrow \P^{N_k}$ be the associated
Kodaira imbedding. Let $Z_k$ be the image of $\iota_\us$. Define
$M_{\al\b}^{(k)}= M_{\al\b}(Z_k)$ as in (\ref{M}). Then

\bea
\label{M1}
 M^{(k)}_{\al\b}\ \ &=& \  \I_{Z_k}
{z_\al\bar z_\b\over |z|^2} \
\o_{FS}^n\ =
\
\I_{X} (s_\al,s_\b)_{h^k} \cdot
\left({h(k)\over h}\right)^k\cdot{\o(k)^n\over
\o^n}\cdot\o^n \ \cr
&=&
\ \d_{\al\b}\ - \ k^{-1}\I_X
(s_\al,s_\b)_{h^k} \cdot s(\o)\ \o^n \ + O(k^{-2})\cr
&=&
\left({1\over n!}\I_X\o^n\right)\cdot{k^n\over N_k}\cdot\d_{\al\b}\ - \ k^{-1}\I
(s_\al,s_\b)_{h^k} \cdot [s(\o)-\bar s]\ \o^n \ + O(k^{-2})\cr
&=&
\left({1\over n!}\I_X\o^n\right){k^n\over N_k}\cdot\d_{\al\b}\ - \ k^{-1}\I
{x_\al\bar x_\b\over |x|^2} \cdot [s(\o)-\hat s]\ \o_{FS}^n \ +
O(k^{-2})\nonumber\\ \cr
\eea
\end{lemma}
where $\hat s\in\R$ is defined by: $\I_X [s(\o)-\hat s]\o^n=0 $.
The proof of Lemma \ref{D05} follows from (\ref{powers of k}).

\v
We return to the proof of Lemma \ref{volest}:
Applying (\ref{M1}) to (\ref{dot zero}) we obtain:

\bea
{1\over n+1}
\label{int} \I_X\dot \phi(0;k)\ \o_{\phi(0;k)}^n  &=&
{1\over k}\cdot{1\over k^n}\cdot \dot E_{\o_{FS}}(\psi_t)
= {1\over k^n}{2\over k}\tr(A_kM^{(k)})\cr\cr
&=& - \
{1\over k^n}{2\over k}
\I\s_\al(s_\al,s_\al)_{h^k}{\l^{(k)}_\al\over k}[s(\o)-\hat s]\ \o^n \ +
O(k^{-2}){1\over k^n}{1\over k}\s_\al |\l_\al|\nonumber
\eea
\be
\ee
where in the first equality we use the fact that $A_k=A_k^*$ and
in the last equality we have made use of the fact that
$A_k$  is traceless.

\v

If we apply Lemma \ref{bound} to equation
(\ref{int}) we obtain

\be
\left|\label{int1} \I_X\dot \phi(0;k)\ \o_{\phi(0;k)}^n\right|\ \leq \
C'{1\over k^n}{1\over k}\sum_\al\I |s_\al|_{h^k}^2\ \o^n\ + \ O(k^{-2}){1\over k^n}{1\over k}CkN_k
\ee
where $C'=C\sup_X|s(\o)-\hat s|$. The first term equals $C'{1\over k}{N_k+1\over k^n}$.
Since ${N_k+1\over k^n}$ is bounded as a function of $k$, the first term is
$O(k^{-1})$. Similarly, the second term is $O(k^{-2})$.

\subsection{The Monge-Amp\`ere equation on a punctured disk}

We now complete the proof of Theorem \ref{main}. As in the proof of Theorem 3,
\cite{PS06}, we can choose a sequence of
positive real numbers $c_k\searrow 0$ in such a way that

\bea\phi(0;k)+c_k\  > \ \phi(0;k+1)+c_{k+1}\eea

Indeed, by the Tian-Yau-Zelditch theorem,
${\rm sup}_X|\phi(0;k)-\phi|\leq C\,k^{-2}$, so that
the sequence
$c_k=2\,C\,\sum_{j\geq k}j^{-2}$ is such a choice. Choose also
$\e_k=k^{-1/2}$ and make the replacement

\bea \label{replacement}
\phi(t;k)\quad\longrightarrow\quad \phi(t;k)+c_k-\e_kt. \eea

Then it
is still true that $\phi(0;k)\to\phi$, and that $\I
\O_{\Phi(k)}^{n+1}=O({1\over k})$. Moreover, the value of $\phi_t$,
as
defined in (\ref{geo}) does not change under this replacement.

\v
Next, we show that $\phi_t$ is continuous at $t=0$
and has the desired initial value.
As in \cite{PS06}, the essential ingredient is a uniform bound for
$|Y\phi(t;k)|$ near the boundary $S^1\times X$, where $Y=\pl_t$.
In fact, differentiating the expression for $\phi(t;k)$ gives
\bea
|\dot\phi(t;k)|
&\leq&
{2\over k}{\sum_{\al=0}^N|\lambda_\al|e^{2t\lambda_\al}|s_\al|^2
\over
\sum_{\al=0}^Ne^{2\lambda_\al}|s_\al|^2}+\e_k
\nonumber\\
&\leq& {2\over k}{\rm sup}_\al |\lambda_\al|+\e_k\,\leq C,
\eea
where in the last step, we made use of the bound $\|A_k\|_{op}\leq C\,k$
provided by Lemma \ref{bound}.
On the boundary $X\times S^1$, the monotonicity of $\phi(t;k)$ guarantees that,
for any pair $k,l$ with $k<l$,
\bea
\phi(t;k)-\phi(t;l)\
>
\
\phi(t;k)-\phi(t;k+1)
\
>\ \delta_k\quad{\rm on}\ X\times S^1,
\eea
where $\delta_k$ is a strictly positive constant independent of $l$.
Since $|\dot\phi(t;m)|$ is uniformly bounded in $m$, it follows that
$\phi(t;k)-\phi(t;l)>{1\over 2}\delta_k$ in a neighborhood $U_k$ of $X\times S^1$
independent of $l$. Thus, we have for any $k$,
\bea
[{\rm sup}_{l\geq k}\phi(t;l)]^*=\phi(t;k)
\eea
in an open neighborhood $U_k$ of $X\times S^1$.
Extend now the original potential $\phi$ on $X$
as a function in a neighborhood of $X\times S^1$, by making it constant along
the flow lines of $Y$. For any $\e>0$, choose $k$ large enough so that
${\rm sup}_X|\phi(0;k)-\phi|<\e$. Then the above estimate
for $\dot\phi(t;k)$ shows that we have
\bea
{\rm sup}_U|\phi(t;k)-\phi|<2\e
\eea
for some neighborhood $U$ of $X\times S^1$ in $X\times D^\times$,
independent of $k$. This implies that
$\phi_t={\rm lim}_{k\to\infty}[{\rm sup}_{l\geq k}\phi(t;k)]^*$
is continuous at $X\times S^1$, and that $\phi_t=\phi$ at $t=0$.

\v

On the other hand,
for fixed $l\geq j>k>0$,  Lemma \ref{phi bound lemma} implies that
\bea
\phi(t;j)-\phi(t;k)
&\leq& C_{k,j}+c_j-c_k+2t f(k)-2tf(j)-(\e_k-\e_j)|t|
\nonumber\\
&\leq&\ [{\rm sup}_{k<m\leq l}C_{k,m}]+1\ -\ {1\over 2}(\e_k-\e_{k+1})|t|.
\eea
for $k$ sufficiently large.
Thus, we can make sure that
\bea
\phi(t;k)>1+\phi(t;j)
\eea
for all $j$ such that $k< j\leq l$ and all $t$ such that
\bea
|t|> 2\,{2+{\rm sup}_{k<m\leq l}C_{k,m}
\over \e_k-\e_{k+1}}\equiv -\log\,r_{k,l}.
\eea
Clearly, we have
$r_{k,l+1}<r_{k,l}$ and, by choosing $C_{k,l}$
in Lemma \ref{phi bound lemma} large enough, we can make sure that
\bea
\lim_{l\to\i}r_{k,l}=\lim_{l\to\i}r_{l,l+1}=0
\eea
for each $k>0$.
Thus, if we set
\bea
\phi(t;k,l)=\sup_{k\leq j\leq l}[\phi(t;j)],
\eea
and let $\O_{k,l}$ be the $(1,1)$ form on $X\times D^\times$
corresponding to $\phi(t;k,l)$ via (\ref{11form}),
we have, for $l>k$,
\be
\label{vol} \I_{X\times D_{r_{k,k+1}}}\O_{k,l}^{n+1}\ \leq \
\I_{X\times D_{r_{k,l}} }\O_{k,l}^{n+1}\ = \
\I_{X\times D_{r_{k,l}} }
\O_{\Phi(k)}^{n+1} \ \leq \
\I_{X\times D^\times }
\O_{\Phi(k)}^{n+1}\ \leq \ {C\over k}
\ee
where $C$ is independent of $k,l$ and
$$ D_{r}\ = \ \{w\in\C: 1>|w|>r\}.
$$
In the middle equality above, we made use of the fact that
the volume integrals depend only
on the values of the currents in a neighborhood of the boundary
of $X\times D_{r_{k,l}}$ (see Lemma 2 in \cite{PS06}).

\v
Now $\phi(t;k,l)$ is an increasing sequence in the index $l$
which converges pointwise, almost everywhere, to
$\xi(t;k)=\sup_{k\leq j}[\phi(t;j)]^*$. Let
$\Xi_k $ be the $(1,1) $ form on $D^\times\times M$
corresponding to $\xi(t;k)$. Then, by the Bedford-Taylor
monotonicity theorem \cite{BT82} applied to the {\it increasing} sequence
$\phi(t;k,l)$ (see also Blocki \cite{B} and Cegrell \cite{Ce}), we
have

\be
\label{vol 2}
\I_{X\times D_{r_{k,k+1}}}\Xi_k^{n+1}\ = \
\lim_{l\to\i}\I_{X\times D_{r_{k,k+1}}}\O_{k,l}^{n+1} \ \leq \ {C\over k}
\ee

Moreover, if $l\geq k$,

\be
\I_{X\times D_{r_{k,k+1}}}\Xi_l^{n+1}\ \leq \
\I_{X\times D_{r_{l,l+1}}}\Xi_l^{n+1}\ \ \leq \
{C\over l}\
\ee

Finally, since $\xi(t;l)$ is monotonically decreasing
to $\phi(t)$ (by definition of $\phi(t)$) we have, using the Bedford-Taylor
monotonicity
theorem again (but this time for decreasing sequences):
\bea
\I_{X\times D_{r_{k,k+1}}}\O_\Phi^{n+1}\ \ =
\lim_{l\to\i}\I_{X\times D_{r_{k,k+1}}}\Xi_l^{n+1}
\  = \ 0
\eea
Since this is
true for all $k$, we obtain
\bea \I_{X\times D^\times}\O_\Phi^{n+1}\ = \ 0.
\eea
Thus $\O_\Phi^{n+1}=0$, and this proves the theorem. Q.E.D.

\section{Proof of Theorem 2}
\setcounter{equation}{0}
\v
In this section we show that if the
expression $N_2(T)^2$ defined below by (\ref{N2}) is strictly positive,
then $\phi_t$ is a non-trivial geodesic, i.e.,
$\phi_t$ is not a constant function of $t$.

\v

Let $N_1+1=\dim(X,L)$ and set $N=N_1$.
Let $\l_0\geq \l_1\geq ...\geq \l_N$ be the diagonal entries of $A_1$ and
let $N_1+1=\dim H^0(X,L)$.
By Lemma \ref{rigid} we may assume that $T$ is imbedded in $\P^N$ and
that the action $\r(\tau)$ is given by the diagonal matrix
whose diagonal entries are given by
$\tau^{\l_0},...,\tau^{\l_N}$.
As usual, we denote by $X_0\sub\P^N$ the central fiber of $T$.
\v

Define $h:\P^N\ra \R$ by
\bea
h(z)={\sum_{\al=0}^N\lambda_\al |z_\al|^2\over\sum_{\al=0}^N|z_\al|^2}.
\eea
\v

We next recall the formula in Donaldson
\cite{D05}:
\bea
{\rm Tr}(A_k^2)=N_2(T)^2\cdot k^{n+2}+O(k^{n+1}),
\eea
where the coefficient $N_2(T)$ is given by
\bea
\label{N2}
N_2(T)^2=\int_{X_0}(h-\hat h)^2\o_{FS}^n,
\eea
and $\hat h$ is determined by
$\int_{X_0}(h-\hat h)\o_{FS}^n=0$. The
test configuration $T$ is trivial if and
only if $N_2(T)=0$.
\v
Now let $\l=\l_N$ so that $\l\leq \l_\al$ for all $\al$.
Denote by $\l_0^{(k)}\geq\l_1^{(k)}\geq\cdots
\geq\l_{N_k}^{(k)}$ the eigenvalues of the endomorphism
$A_k$ (for convenience, the eigenvalues are ordered here in
the opposite order than previously). Set $\l^{(k)}=\l_{N_k}^{(k)}$. Then
$\l^{(k)}=k\l-{{\rm Tr}\,B_k\over N_k}$,
and $\l^{(k)}k^{-1}$ has a limit. Set
\bea
\Lambda={\rm lim}_{k\to\infty}{\l^{(k)}\over k}.
\eea
If $N_2(T)>0$ then (\ref{N2}) implies that the average
absolute eigenvalue $|\l_j^{(k)}|$ has size at least $k$.
On the other hand, since $\Tr A_k=0$, one easily sees
that $|\l^{(k)}|$ has size at least $k$ and thus $\L>0$.
\v
Next, recall that
$$
\phi(t;k)\ =
{1\over k}\log\left(k^{-n}\cdot\sum_{\al=0}^{N_k}\
e^{2t\l^{(k)}_\al}|s^{(k)}_\al|^2_{h_0^k}\right)
$$
Observe that
$$ \I \sum_{\al=0}^{N_k}\ e^{2t\l_\al}|s_\al|^2_{h_0^k} \o_0^n \ \geq \ e^{2|t \l^{(k)}|}
$$
so
\be\label{sup bound}
2|t|{|\l^{(k)}|\over k} + O({1\over k^2})\ \geq \
\sup_{X_t}\phi(t;k)\ \geq\ 2|t| {|\l^{(k)}|\over k}\ - \ n{\log k\over k} -
{1\over k}\log\left(\I_X\o_0^n\right)
\ee

so, letting $k\to\infty$,

\be\label{sup bound}
\sup_{X_t}\phi_t\ = \ 2\,|t|\cdot{|\L|}
\ee

Since $|\L|>0$, this already shows that $\phi_t$ is non-trivial if
$N_2(T)>0$. This establishes Theorem 2. Q.E.D.
\v\v

Under the additional technical assumption (which we expect
can be removed) that for $k_0$ large enough,
$[{\sup}_{k\geq k_0}\phi(t;k)]^*=\sup_{k\geq k_0}\phi(t;k)$ for $|t|>t_{k_0}>>1$,
then the geodesic $\phi_t$ can be shown to be non-trivial
in the stronger sense that it defines a non-trivial ray in ${\cal H}/{\bf R}$.
\v

To show strong non-triviality, we observe that
 $N_2(T)>0$
implies (and is in fact, equivalent to) the following:
There exist $p\in X$ such that $s_\al(p)=0$ for all $\al $ such
that $\l_\al =\l$. Fix such a $p$. Let $\g={\rm inf}\{\l_\al: \l_\al>\l\}$
and $\g^{(k)}={\rm inf}\{\l^{(k)}_\al: \l^{(k)}_\al>\l^{(k)}\}$.
Note that $|\g|<|\l|$.
Again, $\g^{(k)}=k\g-{{\rm Tr}\,B_k\over N_k}$, and
${\g^{(k)}\over k}$ has a limit as $k\to\infty$,
\bea
\Gamma={\lim}_{k\to\infty} {\g^{(k)}\over k}
\eea
satisfying $|\Gamma|<|\Lambda|$. Now, at the point $p$, we have for all $k$,
\be
\phi(t;k)(p)\ \leq \ 2|t|{|\g^{(k)}|\over k} \ + \ O({1\over k}).
\ee
Fix $\epsilon>0$ so small that $|\Gamma|+2\e<|\Lambda|$. Then there exists $k_0$ so that
\bea
\phi(t;k)(p)\ \leq\ 2(|\Gamma|+\e)|t|
\eea
for all $|t|>1$ and all $k\geq k_0$. We have then, for $|t|$ sufficiently large,
\bea
\phi_t(p)
\leq [{\rm sup}_{k\geq k_0}\phi(t;k)]^*(p)\ = \
{\rm sup}_{k\geq k_0}\phi(t;k)(p)\
\leq \ 2(|\Gamma|+\e)|t|.
\eea
In view of (\ref{sup bound}),
this shows $\lim_{t\to -\i}\osc_{X_t} \phi_t = \i$ where
 $\osc_{X_t} \phi_t = \sup_{X_t} \phi_t - \inf_{X_t}\phi_t$.
Thus $\phi_t$ is strongly non-trivial.

\section{Proof of Theorem 3}
\setcounter{equation}{0}

The formula in Lemma 8.8 of Tian \cite{T97}
implies that

\be \lim_{t\to -\i}\dot \n_k\ = \ F_{CM}(T)
\ee
where $F_{CM}(T)$ is the CM-Futaki invariant
(see \cite{T97} for the precise definition).
\v
On the other hand, the recent work of Paul-Tian
\cite{PT06} shows that $F_{CM}(T)=F(T)$
under the hypothesis of Theorem 3.


\newpage

\begin{thebibliography}{99}

\bibitem{AT} Arezzo, C., and Tian, G.,
``Infinite geodesic rays in the space of K\"ahler potentials'',
Ann. Sci. Norm. Sup. Pisa Sci. (5) 2 (2003) 617-630,
arXiv : math.DG / 0210389.

\bibitem{BT76} Bedford, E. and B.A. Taylor,
``The Dirichlet problem for a complex Monge-Amp\`ere equation'',
Invent. Math. {\bf 37} (1976), 1-44.

\bibitem{BT82} Bedford, E. and B.A. Taylor,
``A new capacity for plurisubharmonic functions'', Acta Math. {\bf 149}
(1982), 1-40.


\bibitem{Be05} Berndtsson, B.,
``Curvature of vector bundles associated to holomorphic fibrations''
preprint CV/0511225

\bibitem{Be06} Berndtsson, B., ``Positivity properties of direct
image bundles'', preprint.


\bibitem{B} Blocki, Z.,
``The complex Monge-Amp\`ere operator and pluripotential
theory'', lecture notes available from the author's website.

\bibitem{C} Catlin, D.,
``The Bergman kernel and a theorem of Tian'',
{\it Analysis and geometry in several complex variables} (Katata,
1997), 1-23,
Trends Math., Birkh\"auser Boston, Boston, MA, 1999.

\bibitem{Ce} Cegrell, U.,
{\it ``Capacities in complex analysis"},
Aspects of Math. E14, Vieweg, 1988.

\bibitem{Ch} Chen, X.X.,
``The space of K\"ahler metrics'', J. Differential Geom.
{\bf 56} (2000), 189-234.

\bibitem{CT} Chen, X.X. and G. Tian,
``Geometry of K\"ahler metrics and foliations by discs'',
arXiv: math.DG / 0409433.

\bibitem{ChIII} Chen, X.X.,
``Space of K\"ahler metrics III - On the lower bound of the Calabi
 energy and geodesic distance", arXiv: math.DG/0606228




\bibitem{D99} Donaldson, S.K.,
``Symmetric spaces, K\"ahler geometry, and Hamiltonian
dynamics'', Amer. Math. Soc. Transl. {\bf 196} (1999) 13-33.

\bibitem{D01} Donaldson, S.K.,
``Scalar curvature and projective imbeddings I'',
J. Differential Geom. {\bf 59} (2001) 479-522.

\bibitem{D02} Donaldson, S.K.,
``Scalar curvature and stability of toric varieties'',
J. Differential Geom. {\bf 62} (2002), 289-349.

\bibitem{D04} Donaldson, S.K.,
``Scalar curvature and projective imbeddings II'',
arXiv: math.DG / 0407534.

\bibitem{D05} Donaldson, S.K.,
 ``Lower bounds on the Calabi functional'',
J. Differential Geom.  {\bf 70}  (2005),  453-472.

\bibitem{K} Kolodziej, S.,
 ``The complex Monge-Amp\`ere equation",
Acta Math.  {\bf 180}  (1998),  no. 1, 69-117.







\bibitem{Lu}  Lu, Z., ``On the lower order terms of the
asymptotic expansion of Tian-Yau-Zelditch'',
Amer. J. Math.  {\bf 122} (2000),   235-273.


\bibitem{M87} Mabuchi, T.,
``Some symplectic geometry on compact K\"ahler manifolds'',
Osaka J. Math. {\bf 24} (1987) 227-252.

\bibitem{Mumford} Mumford, D.,
`` Stability of projective varieties'',
Enseignement Math. {\bf 23}  (1977),  39-110

\bibitem{P00} Paul, S.,
``Geometric analysis of Chow Mumford stability'',
 Adv. Math. {\bf 182} (2004), 333-356.

\bibitem{PT06} Paul, S. and G. Tian,
``CM stability and the generalized Futaki invariant'',
arXiv: math: AG/0605278.

\bibitem{PT06II} Paul, S. and G. Tian,
``CM stability and the generalized Futaki invariant II'',
arXiv: math.DG/0606505

\bibitem{PS02} Phong, D.H. and J. Sturm,
``Stability, energy functionals, and K\"ahler-Einstein
metrics'', Comm. Anal. Geometry {\bf 11} (2003)
563-597, arXiv: math.DG / 0203254.

\bibitem{PS02a} Phong, D.H. and J. Sturm,
``Scalar curvature, moment maps, and the Deligne pairing'',
Amer. J. Math. {\bf 126} (2004) 693-712,
arXiv: math.DG /
0209098.

\bibitem{PS03} Phong, D.H. and J. Sturm,
``The Futaki invariant and the Mabuchi energy of a complete
intersection'', Comm. Anal. Geometry {\bf 12} (2004)
321-343, arXiv: math.DG / 0312529.

\bibitem{PS06} Phong, D.H. and J. Sturm,
``The Monge-Amp\`ere operator and geodesics in the space of K\"ahler potentials'',
Invent. Math. 2006 (to appear),
arXiv: math.DG / 0504157.





\bibitem{RT} Ross, J. and Thomas, R.,
``A study of the Hilbert-Mumford criterion for the stability of projective varieties''
arXiv : math.AG /0412519

\bibitem{S92} Semmes, S.,
``Complex Monge-Amp\`ere and symplectic manifolds'',
Amer. J. Math. {\bf 114} (1992) 495-550.

\bibitem{T90} Tian, G.,
`` On a set of polarized K\"ahler metrics on algebraic
manifolds'', J. Diff. Geom. {\bf 32} (1990) 99-130.


\bibitem{T97} Tian, G.,
``K\"ahler-Einstein metrics with positive scalar curvature'',
Invent. Math. {\bf 130} (1997) 1-37.


\bibitem{Y78} Yau, S.T.,
``On the Ricci curvature of a compact K\"ahler manifold
and the complex Monge-Amp\`ere equation I'',
Comm. Pure Appl. Math. {\bf 31} (1978) 339-411.

\bibitem{Y87} Yau, S.T.,
`` Nonlinear analysis in geometry'',
Enseign. Math. (2)  {\bf 33}  (1987),  no. 1-2, 109--158.


\bibitem{Y} Yau, S.T.,
``Open problems in geometry'',
Proc. Symp. Pure Math. {\bf 54}, AMS Providence, RI
(1993) 1-28.

\bibitem{Z} Zelditch, S.,
``The Szeg\"o kernel and a theorem of Tian'',
Int. Math. Res. Notices {\bf 6} (1998) 317-331.

\bibitem{Z96} Zhang, S.,
``Heights and reductions of semi-stable varieties'',
Compositio Math. {\bf 104} (1996) 77-105.


\end{thebibliography}
\enddocument